\documentclass[12pt]{amsart}
\usepackage{graphicx, amsmath, amsfonts, amssymb, amsthm,geometry,xcolor,enumerate}
\usepackage[hidelinks]{hyperref}
\usepackage{fullpage}

\newtheorem{problem}{Problem}

\newtheorem{theorem}{Theorem}[section]
\newtheorem{corollary}{Corollary}[section]

\newtheorem*{RH}{Riemann Hypothesis (RH)}

\newtheorem*{density}{Density Hypothesis}
\newtheorem*{littlewood}{Littlewood's Lemma}
\newcommand{\re}{\textnormal{Re}}
\newcommand{\im}{\textnormal{Im}}
\newcommand{\Li}{\textnormal{Li}}
\newcommand{\Arg}{\textnormal{arg}\,}

\title{A decades-long breakthrough in\\ zero-density estimates and primes in short intervals\\}

\author[C.~L. Turnage-Butterbaugh]{Caroline~L.~Turnage-Butterbaugh}

\dedicatory{To the memory of Fredric T. Howard}

\begin{document}

\begin{abstract}
The Riemann Hypothesis (RH) asserts that every nontrivial zero of the Riemann zeta-function has real part equal to $1/2$. A zero-density theorem provides evidence towards RH by bounding the number of zeros of the zeta-function with real part greater than $1/2$. In 2024, Larry Guth and James Maynard announced a new zero-density theorem which, for a key location in the critical strip, strengthens previous work of Ingham and is the first such improvement in over 80 years. This expository paper places this remarkable achievement in the context of the rich history of zero-density theorems and explores its implications on the distribution of primes in short intervals.
\end{abstract}
\maketitle
\section{Introduction}
In 1896, Hadamard \cite{Hadamard96} and de la Vallée Poussin \cite{VP1896} independently proved the Prime Number Theorem, which describes the asymptotic distribution of the primes among the integers. Both proofs stemmed from the foundational work introduced decades earlier by Riemann \cite{Rieman1859} in which it was demonstrated that the distribution of the primes is fundamentally related to the complex zeros of the now eponymous Riemann zeta-function. The manuscript, Riemann's only paper in number theory and six pages in length, is also the birthplace of one of the most important open questions in mathematics: the Riemann Hypothesis, which asserts that the complex zeros of the Riemann zeta-function all lie on the \emph{critical line} $\re(s)=1/2$. 

Evidence in favor of the Riemann Hypothesis comes in many forms. One type of result concerns the proportion of zeros lying on the critical line.  To date it is known \cite{Pratt20} that at least five-twelfths of the complex zeros of the Riemann zeta-function satisfy the Riemann Hypothesis. A second type of evidence centers on showing that there are zero-free regions in the \emph{critical strip}, where $0\le \re(s) \le1$. The strongest zero-free region known, due to Vinogradav \cite{Vinogradov58} and Korbov \cite{Korobov58}, is discussed at the end of Section \ref{primes}. The striking recent work of Guth and Maynard (announced in ~\cite{GuthMaynard24} and now published in \cite{GuthMaynard2026}) gives a breakthrough on another type of evidence in favor of the Riemann Hypothesis, called a zero-density estimate, which bounds the size of the (conjecturally empty) set of zeros which lie more than a fixed distance from the critical line. Guth and Maynard proved a zero-density theorem with a substantially improved bound that, for a key location in the critical strip, had not seen improvement since the work of Ingham \cite{Ingham40} in 1940. Zero-density estimates are also a key tool in the study of the distribution of the primes in short intervals. Guth and Maynard's work gives the first substantial improvement for this application since the work of Huxley \cite{Huxley72} in 1972. 

In this paper we survey the rich history of zero-density estimates and primes in short intervals. We begin with a brief sketch of the proof of the Prime Number Theorem to make explicit the connection between the zeros of the Riemann zeta-function and the distribution of the primes. We then demonstrate the necessity of a zero-density estimate for the study of the distribution of primes in short intervals, and examine how Dirichlet polynomials play a central role in obtaining these estimates. Finally, we discuss Guth and Maynard's new ideas and insights.\\

\noindent\textbf{Notation.}\,Throughout this exposition we follow Riemann, setting $s=\sigma+it$, where $\sigma,t$ are real. We also make use of the following notation. For $g(x)>0$, $f(x)=O(g(x))$ and $f(x) \ll g(x)$ both mean that there exists a constant $M>0$ such that $|f(x)|\le Mg(x)$ for all $x$ sufficiently large. We write $f(x)\sim g(x)$ to mean  $\lim_{x\to \infty}f(x)/g(x)=1$ and $f(x)=o(g(x))$ to mean  $\lim_{x\to \infty}f(x)/g(x)=0$. An integral written as $\int_{(c)}$ means $\lim_{T\to \infty}\int_{c-iT}^{c+iT}$. Finally, we frequently make use of $\varepsilon$ to denote a small, positive quantity, which may change from line to line.

\section{The distribution of the primes and the zeros of the Riemann zeta-function}\label{primes} Let $x>0$, and consider the prime counting function
\[
\pi(x) := \#\left\{p \text{ prime} \,\big|\, p \le x\right\}.
\]
Since there are infinitely many primes,  it follows that $\pi(x)\to \infty$ as $x\to \infty$. But \emph{how} does $\pi(x)$ grow? Gauss (1792) conjectured\footnote{Legendre \cite{Legendre1798} conjectured that $\pi(x) \sim x/(\log x + 1.08366)$ as $x\to \infty$. Both Gauss' and Legendre's conjectures imply that  $\pi(x) \sim x/\log x$ as $x\to \infty$.}
\begin{equation}\label{eqn: pnt}
\pi(x) \sim \Li(x) := \int_{2}^{x}\frac{du}{\log u}, \quad x \to \infty.
\end{equation}
Repeated integration by parts shows that $\Li(x)\sim x/\log x$ as $x\to \infty$. In 1896, Hadamard and de la Vallée Poussin independently proved that the asymptotic in \eqref{eqn: pnt} holds, and we now refer to it as the \emph{Prime Number Theorem}. Their proofs were made possible by Riemann's groundbreaking 1859 manuscript which introduced new analytic tools and ideas to study the primes via properties of the Riemann zeta-function. This function is defined by the Dirichlet series
\begin{equation}\label{def: zeta}
\zeta(s) := \sum_{n=1}^{\infty}\frac{1}{n^s},\qquad \re(s)>1
\end{equation}
and has meromorphic continuation to the entire complex plane with one singularity: a simple pole at $s=1$. As we will see below, at the heart of the proof of the Prime Number Theorem is the fact that $\zeta(s)$ does not vanish when $\re(s)=1$, and stronger information about the location of the nontrivial zeros of $\zeta(s)$ yields stronger forms of the Prime Number Theorem in which error terms are specified. Before we sketch how this works, we record a few important properties of $\zeta(s)$.

\subsection{Properties of \texorpdfstring{$\zeta(s)$}{zeta(s)}} One fundamental property of $\zeta(s)$ is its functional equation, which states that for all complex $s$, 
\begin{equation*}
\xi(s)=\xi(1-s),
\end{equation*}
where
\[
\xi(s):=\frac{1}{2}s(s-1)\pi^{-s/2}\Gamma(s/2)\zeta(s),
\]
and $\Gamma(s)$ denotes the Gamma function. The functional equation shows that $\zeta(s)$ is symmetric about the critical line. Since $\xi(s)$ is entire, it follows that $\zeta(s)$ has zeros, called \emph{trivial zeros}, at the poles of $\Gamma(s)$, i.e. at $s=-2n$ for $n\in \mathbb{N}$. We also note that $\zeta(s)$ has infinitely many \emph{nontrivial zeros} in the critical strip.\footnote{We will often refer the nontrivial zeros of $\zeta(s)$ as `zeros.'} In fact, the \emph{Riemann--von~Mangoldt Formula}, which was conjectured by Riemann~~\cite{Rieman1859} and proved by von~Mangoldt ~\cite{vonMangoldt} in 1905, provides a surprisingly precise formula for $N(T)$, the number of nontrivial zeros in the rectangle ${0 \le \re(s) \le 1}$ and  $0\le \im(s) \le T$. It states 
\begin{equation}\label{eqn: NT}
N(T) = \frac{T}{2\pi}\log \frac{T}{2\pi}-\frac{T}{2\pi}+O(\log T), \quad T\to \infty.
\end{equation}
How does one prove such a formula? You may recall the argument principle from complex analysis, which connects the difference between the number of zeros and poles of a meromorphic function in a region to a contour integral of the function's logarithmic derivative. See Chapter IX, Section 9.3 of ~\cite{Titchmarsh} for the details on applying the argument principle to count the nontrivial zeros of $\zeta(s)$ in the rectangle.

\subsection{A brief sketch of counting primes}
How can it be that the location of the nontrivial zeros of $\zeta(s)$ has a fundamental impact on the distribution of the primes? To begin to answer this question, we introduce the Chebyshev $\psi$-function, a modification of $\pi(x)$ which counts prime powers $p^k \le x $ with the weight $\log p$. Let
\[
\psi(x) := \sum_{n\le x}\Lambda(n),  
\]
where $\Lambda(n)$ is the von Mangoldt function given by
\[
\Lambda(n):=
    \begin{cases}
    \log p &\text{if } n = p^k,\, k \in \mathbb{N},\\
    0 & \text{else}.
    \end{cases}
\]
We can  control the contribution from the prime powers with $k\ge 2$ as follows: $p^k \le x$ implies $k  \le \log_2 x$, and therefore
\[
\sum_{\substack{p^k\le x\\k\ge 2}}\log p = \sum_{2 \le k \le \log_2x}\left( \sum_{p \le x^{1/k}}\log p\right)\le 
\sum_{2 \le k \le \log_2x}x^{1/k}\log x
\le x^{1/2}\log^2 x = o(x), \quad x\to \infty.
\]
Thus for large $x$  we should (roughly) have
\[
\psi(x) = \sum_{n\le x} \Lambda(n) \approx \sum_{p\le x}\log p \approx \pi(x)\log x.
\]
Making this argument precise shows that the statement $\pi(x) \sim x / \log x$ as $x \to \infty$ is equivalent to the statement that $\psi(x) \sim x$ as $x\to \infty$. (See, e.g., ~\cite[Theorem 4.4]{Apostol76}.)

We now sketch an argument of Riemann which shows that the asymptotic behavior of $\psi(x)$ explicitly depends on the nontrivial zeros of $\zeta(s)$. By the Fundamental Theorem of Arithmetic, the Dirichlet series \eqref{def: zeta} defining $\zeta(s)$ may be expressed as the Euler product
\begin{equation}\label{def: eulerprod}
\zeta(s)=\prod_{p \text{ prime}}\left(1-\frac{1}{p^s}\right)^{-1}, \quad \re(s) >1.
\end{equation}
The Euler product reveals that $\zeta(s)$ has no zeros for which $\re(s)>1$. Taking the logarithmic derivative of \eqref{def: eulerprod} and expressing $1-p^{-s}$ as a geometric series yields
\begin{equation}\label{eqn: logderivzeta}
-\frac{\zeta'(s)}{\zeta(s)} = \sum_{p} \frac{\log p}{p^s(1-p^{-s})}= \sum_{p}\log p \left(\sum_{m=1}^{\infty}\frac{1}{p^{ms}} \right) = \sum_{\substack{p \\ m \ge 1}}\frac{\log p}{p^{ms}} = \sum_{n=1}^{\infty}\frac{\Lambda(n)}{n^s}.
\end{equation}
We may extract $\psi(x)$ from \eqref{eqn: logderivzeta} using Perron's formula, or, said another way, via the discontinuous integral
\begin{equation}\label{eqn: perronintegral}
\frac{1}{2\pi i}\int_{(c)}\,\frac{y^s}{s}\,ds = \begin{cases}
1 \quad \text{if } y>1,\\[.1in]
\displaystyle\frac{1}{2} \quad \text{if } y=1,\\[.1in]
0 \quad \text{if } 0<y<1,\\
\end{cases}
\end{equation}
where $c>1$. If $x$ is not an integer,\footnote{Of course, $x\in \mathbb{N}$ is a possibility! In this case, one must adjust the last term in the sum defining $\psi(x)$, replacing $\Lambda(x)$ with $\frac{1}{2}\Lambda(x)$.} then by \eqref{eqn: logderivzeta} and \eqref{eqn: perronintegral}  
\[
\psi(x) = \sum_{n\le x}\Lambda(n) =\sum_{n=1}^{\infty}\Lambda(n)\left(\frac{1}{2\pi i}\int_{(c)}\,\frac{(x/n)^s}{s}\,ds\right) =\frac{1}{2\pi i}\int_{(c)}\,-\frac{\zeta'(s)}{\zeta(s)}\frac{x^s}{s}\,ds, 
\]
that is
\[
\psi(x) =-\frac{1}{2\pi i}\int_{(c)}\,\frac{\zeta'(s)}{\zeta(s)}\frac{x^s}{s}\,ds.
\]
Thus, we can study $\psi(x)$ via the integral above. In particular, we may use the calculus of residues to ``shift the contour" all the way to infinity on the left, thereby expressing $\psi(x)$ in terms of the poles of the integrand. The residue from the simple pole at $s=1$ contributes $x$. The pole at $s=0$ contributes $-\zeta'(0)/\zeta(0)$ which is equal to $-\log(2\pi)$. (See (8) and (10) of \cite[Chapter 12]{Davenport}.) The trivial zeros at the negative even integers contribute
\[
-\sum_{n=1}^{\infty}\frac{x^{-2n}}{-2n} = -\frac{1}{2}\log\left( 1-x^{-2}\right),
\]
and each nontrivial zero $\rho$ of $\zeta(s)$ contributes $x^\rho \rho^{-1}$. In all, we deduce the \emph{Explicit Formula}:  for $x>1$,  
\begin{equation}\label{eqn: explicitformula0}
 \psi(x) = x - \lim_{T\to \infty}\sum_{\substack{\rho\\|\gamma|\le T}}\frac{x^{\rho}}{\rho} -\log{2\pi} -\frac{1}{2}\log(1-x^{-2}).
 \end{equation}
We have reached the crux of the argument: since $$|x^\rho|\le x^{\re(\rho)},$$ the contribution from the nontrivial zeros in \eqref{eqn: explicitformula0} is small provided $\re(\rho)\ne1$. (See \cite[Theorem 30]{Ingham} for details on handling  $\sum|\rho|^{-1}$.)  Hadamard and de la Vallée Poussin proved this nontrivial fact (and slightly more) to conclude $\psi(x)\sim x$ as $x\to \infty$. To obtain an error term, one must prove a zero-free region of $\zeta(s)$ in the critical strip which, in turn, gives a stronger upper bound on $\re(\rho)$. Since the zeros are symmetric about the critical line, the best possible scenario is that every zero of $\zeta(s)$ is \emph{on} the critical line. This hypothetical is, of course, the famous Riemann Hypothesis.

\begin{RH}
Every nontrivial zero of $\zeta(s)$ has real part equal to $1/2$.
\end{RH}
If the Riemann Hypothesis is true, then 
\begin{equation}\label{eqn: pntRH}
\pi(x) = \Li(x) + O(x^{1/2}\log x), \quad x\to \infty.
\end{equation}
Let $s=\sigma+it$. We do not currently know of a single value $\delta>0$ for which $\zeta(s)$ has no zeros in the region $\sigma \ge 1-\delta$.  All known zero-free regions of $\zeta(s)$ are of the form 
\[
\sigma \ge 1-c\,\eta(t),
\]
where $c>0$ and, as $t\to \infty$,  $\eta(t)\to 0$ as illustrated below. 
\begin{figure}[ht]
    \centering
    \includegraphics[width=0.4\linewidth]{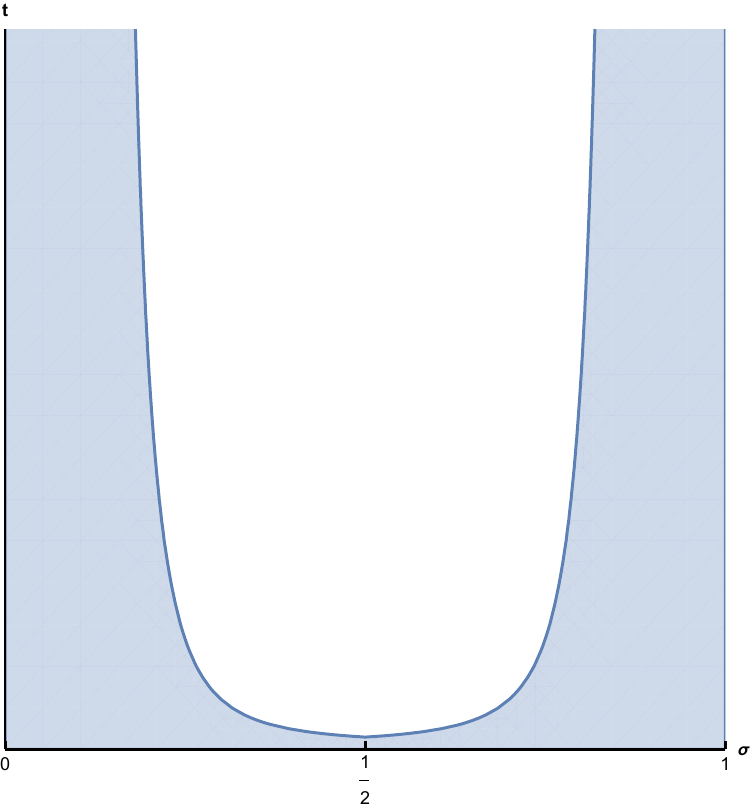}
    \caption{The general shape of known zero-free regions (in blue) for $\zeta(s)$ in the upper half-plane of the critical strip.}
\end{figure}

In 1899, de la Vallée Poussin \cite{VP1899} proved that for sufficiently large $t$ there exists an absolute constant $c>0$ such that $\zeta(s)$ has no zeros in the region 
\begin{equation}\label{eqn: zerofreeclassical}
\sigma \ge 1 - c/\log t.
\end{equation}
The best unconditional zero-free region is due to Vinogradov ~\cite{Vinogradov58} and Korobov ~\cite{Korobov58}. They proved that for sufficiently large $t$, there exists an absolute constant $c>0$ such that $\zeta(s)$ has no zeros in the region 
\begin{equation}\label{eqn: zerofreeVK}
\sigma \geq 1-c(\log t)^{-2/3}(\log\log t)^{-1/3}.
\end{equation}
This, in turn, implies
\[
\pi(x) = \Li(x) + O\left(x\exp\left(-c(\log x)^{3/5}(\log\log x)^{-1/5} \right) \right).
\]

\section{The distribution of the primes in short intervals} Let $x, y$ be large, where $y=y(x)$ is smaller than $x$. How many primes are there in the interval $(x,x+y]$? Since there are roughly $y$ integers in the interval and each has a $1/\log x$ probability of being prime by \eqref{eqn: pnt}, one might expect that there are roughly $y/\log x$ primes in the interval. In other words, for suitable choices of $x$ and $y$, we expect 
\[
\pi(x+y)-\pi(x) \approx  \frac{y}{\log x},
\]
or, equivalently, 
\[
\psi(x+y)-\psi(x) \approx  y.
\]
More precisely, we say that the Prime Number Theorem holds in the short interval $(x, x+y]$ if 
\begin{equation}\label{eqn: pntshort}
\psi(x+y) - \psi(x) = \sum_{x < n \le x+y}\Lambda(n) \sim y, \qquad x \to \infty.
\end{equation}
For the purposes of this exposition, let $y=x^{\theta}$, where $0 < \theta \le 1$, and note that the case $\theta=1$ corresponds to the Prime Number Theorem. We thus consider the following problem. 

\begin{problem}\label{problemA}
Let $y=x^\theta$ and $0 < \theta < 1$. Determine the least exponent $\theta$ for which \eqref{eqn: pntshort} holds for sufficiently large $x$.
\end{problem}
If we assume the Riemann Hypothesis, then \eqref{eqn: pntRH} yields 

\[
\pi(x+x^\theta)-\pi(x) = \int_{x}^{x+x^\theta}\frac{du}{\log u} + O\left(x^{1/2}\log x\right). 
\]
Note that the error term is smaller than the main term provided $\theta > 1/2 + \varepsilon$, where $\varepsilon >0$. How small can we take $\theta$ unconditionally? 

As a first pass, let's proceed as in the proof of the Prime Number Theorem. Let $\rho=\beta+i\gamma$ denote a nontrivial zero of $\zeta(s)$, where $\beta, \gamma \in \mathbb{R}$. By the truncated form of \eqref{eqn: explicitformula0} (see \cite[Chapter 17]{Davenport}) we have
\[
\psi(x+y)-\psi(x) = y-\sum_{\substack{\rho\\|\gamma|\le T}}\frac{(x+y)^\rho-x^\rho}{\rho}+O\left(\frac{x^{1+\varepsilon}}{T} \right), \qquad T\to\infty.
\]
As expected, we must control the contribution from the nontrivial zeros $\rho$ of $\zeta(s)$. Following the set-up of an argument by Ingham ~\cite{Ingham37},  applications of the Fundamental Theorem of Calculus and the triangle inequality yield
\[
\left|\frac{(x+y)^\rho-x^\rho}{\rho}\right|=\left|\int_{x}^{x+y}u^{\rho-1}\,du\right| \le \int_{x}^{x+y}u^{\beta-1}\,du\le yx^{\beta-1}.
\]
Thus 
\begin{equation}\label{eqn: pntshortmid}
\frac{\psi(x+y)-\psi(x)}{y} = 1+O\left(\sum_{\substack{\rho\\|\gamma|\le T}}x^{\beta-1}\right)+O\left(\frac{x^{1+\varepsilon}}{yT} \right), \qquad T\to \infty.
\end{equation}

In order for the expected asymptotic to hold, we need to show that the two error terms on the right-hand side of \eqref{eqn: pntshortmid} are both $o(1)$. Recall that $y=x^\theta$, where $0<\theta<1$, and set $T=x^\alpha$, where $0 < \alpha < 1$. Note that the second error term is $o(1)$ provided $\theta$ and $\alpha$  satisfy
\begin{equation}\label{eqn: sum1}
1-\alpha+\varepsilon < \theta <1.
\end{equation}
The first error term on the right-hand side of \eqref{eqn: pntshortmid} depends on the zeros of $\zeta(s)$. Let us see what happens when, as in the proof of the Prime Number Theorem, we apply a zero-free region to bound the sum.

\subsection{Applying a zero-free region} We illustrate the argument using the classical zero-free region \eqref{eqn: zerofreeclassical}, which guarantees that any nontrivial zero $\rho = \beta + i\gamma$ with $|\gamma|\le T$ and $T$ large satisfies
$\beta < 1 - c\,\eta(T)$, where $\eta(T)=1/\log T$ and $c>0$ is a constant. Applying \eqref{eqn: NT} with $T=x^\alpha$, we find  
\[
\sum_{\substack{\rho\\|\gamma|\le T}}x^{\beta-1} \ll x^{-c\,\eta(T)}T\log T \ll x^{-c\,\eta(x^\alpha)+\alpha+\varepsilon}, \quad x\to \infty.
\]
We have reached an impasse: the classical zero-free region is not strong enough to ensure that the exponent $-c\,\eta(x^\alpha)+\alpha +\varepsilon<0$ for large $x$ because $\eta(x^\alpha)=(\alpha\log x)^{-1}\to 0$ as $x\to \infty$. Note that the Vinogradov-Korobov zero-free region \eqref{eqn: zerofreeVK} shares this feature. In short, no matter the choice of $\alpha$, any application of a currently available zero-free region will fail to show $\sum_{|\gamma|\le T}x^{\beta-1} =o(1)$ as $x, T\to \infty$. On the other hand, if RH is true, then we can take $\alpha < 1/2 - \varepsilon$ and, via \eqref{eqn: sum1}, conclude (again) that RH implies $\theta > 1/2 + \varepsilon$ in Problem \ref{problemA}. 

The main takeaway from this exercise is the observation that the zeros close to the critical line contribute much less to the sum than those close to the line $\re(s)=1$; the above argument failed because we treated each of the roughly $ T\log T$ zeros as if they are in this latter category. In 1930 Hoheisel ~\cite{Hoheisel30} circumvented this issue by leveraging the fact that there are very few zeros of $\zeta(s)$ close to the line $\re(s)=1$.\footnote{Of course, if RH is true, this hypothetical set of ``bad" zeros is empty!} This type of result is made quantitatively precise via what is known as a \emph{zero-density estimate}. 

\subsection{Applying a zero-free region and a zero-density estimate}\label{sec: zfzd}  Let $N(\sigma,T)$ denote the number of zeros $\zeta(s)$ for which $\beta \ge \sigma$ and $0<\gamma \le T$. By \eqref{def: eulerprod}, $N(1,T)=0$, and for any $0\le \sigma<1$, \eqref{eqn: NT} yields $N(\sigma, T) \le N(T)\sim (T/2\pi)\log T$. Under RH, $N(\sigma,T)=0$ for $1/2<\sigma<1$. Therefore we aim to bound $N(\sigma,T)$ in this range.

An upper bound for $N(\sigma,T)$ of the form 
\[
N\left(\sigma,T\right)\ll T^{A(1-\sigma)}\log^BT, \quad T\to \infty,
\]
where $A,B\ge  0$ is called a \emph{zero-density estimate}. We will discuss zero-density estimates more thoroughly in the next section,  but for the moment let's sketch how to use such a result to prove the Prime Number Theorem in a short interval. 

As discussed in ~\cite{Ingham37} (see, also, Section 10.5 of ~\cite{IK04}), Hoheisel's work required two key inserts. Firstly, he made use of a precursor to \eqref{eqn: zerofreeVK} due to Littlewood ~\cite{Littlewood22},  which states that for sufficiently large $t$, $\zeta(s)$ is free of zeros in the region
\begin{equation}\label{eqn: littlewood}
\sigma \ge 1 - c\frac{\log\log t}{\log t},
\end{equation}
where the constant $c>0$ is very small. Secondly, he made use of Carlson's zero-density theorem ~\cite{Carl21}, which for $1/2< \sigma < 1$ implies
\begin{equation}\label{eqn: carlson}
N\left(\sigma,T\right) \ll T^{\,4(1-\sigma)}\log^6 T, \quad T\to\infty.
\end{equation}
With these results in hand, let us return to \eqref{eqn: pntshortmid}, recalling that we aim to show that the sum over zeros there is $o(1)$. Towards this end, note that 
\[
\sum_{\substack{\rho\\|\gamma|\le T}}(x^{\beta-1}-x^{-1}) = \sum_{\substack{\rho\\|\gamma|\le T}}\int_{0}^{\beta}x^{\sigma-1}\log x\,d\sigma
=\int_{0}^{1}\sum_{\substack{\rho\\\substack{|\gamma|\le T \\ \beta \ge\sigma}}}x^{\sigma-1}\log x\,d\sigma,
\]
where the first equality follows from the Fundamental Theorem of Calculus. Thus by \eqref{eqn: NT},
\begin{align*}
\sum_{\substack{\rho\\|\gamma|\le T}}x^{\beta-1} &\le 2x^{-1}N(T)+2\int_{0}^{1/2}N(T)x^{\sigma-1}\log x\,d\sigma+2\int_{1/2}^{1}N(\sigma,T)x^{\sigma-1}\log x\,d\sigma\\
&\ll x^{-1}T\log T +x^{-1/2}T\log T+ \int_{1/2}^{1-c\,\eta(T)}N(\sigma,T)x^{\sigma-1}\log x\,d\sigma,
\end{align*}
where $c\,\eta(T)=c\,\log\log T / \log T$ as in \eqref{eqn: littlewood}. Applying the zero-density estimate \eqref{eqn: carlson} and making the change of variable $\sigma \mapsto 1-\sigma$, we find
\begin{align}
\sum_{\substack{\rho\\|\gamma|\le T}}x^{\beta-1} &\ll x^{-1/2}T\log T + (\log x)\log^6T\int_{c\,\eta(T)}^{1/2}\left(\frac{T^4}{x}\right)^{\sigma}\,d\sigma\notag \\
&\ll x^{-1/2}T\log T+ \log^6T\left(\frac{T^4}{x}\right)^{c\,\eta(T)}.\label{eqn: step}
\end{align}
As in our first attempt to address Problem \ref{problemA}, we aim to pick $\theta$ and $\alpha$ such that \eqref{eqn: sum1} holds. Set $T=x^{\alpha}$, where $0<\alpha <1/4$. Then by \eqref{eqn: step}, 
\[
\sum_{\substack{\rho\\|\gamma|\le T}}x^{\beta-1} \ll x^{\alpha-1/2}\log x + x^{(4\alpha-1)c\,\eta(x^\alpha)}\log^6x.
\]
Moreover, $x^{\alpha-1/2}\log x \ll (\log x)^{-1}$ and 
\[
x^{(4\alpha-1)c\,\eta(x^\alpha)} = \exp\left(\frac{c(4\alpha-1)}{\alpha}\log(\alpha\log x)\right) \ll (\log x)^{\frac{c(4\alpha-1)}{\alpha}}.
\]
Collecting these estimates and setting $\alpha = 1/(4+7c^{-1})<1/4$, we find
\[
\sum_{\substack{\rho\\|\gamma|\le T}}x^{\beta-1}  \ll \frac{1}{\log x}.
\]
Note that by \eqref{eqn: sum1}, $\theta$ depends on $\alpha$ and thus on $c$. Using the zero-free region \eqref{eqn: littlewood}, for which we recall $c>0$ is very, very small, Hoheisel was able to take $\theta = 32999/33000$. If we instead use the Vinogradov-Korobov zero-free region \eqref{eqn: zerofreeVK}, this argument yields $\theta = 3/4 + \varepsilon$.

\section{Zero-density theorems for the Riemann zeta-function}
We noted in Section \ref{sec: zfzd} that the Riemann Hypothesis implies $N(\sigma,T)=0$ for $1/2<\sigma<1$. In 1937, Ingham ~\cite[Theorem 3]{Ingham37} proved that the Lindelöf Hypothesis\footnote{The Lindelöf Hypothesis conjectures that for every $\varepsilon>0$, $\zeta(1/2+it)\ll t^\varepsilon$ as $t\to \infty$.} implies 
\[
N(\sigma,T)\ll T^{2(1-\sigma)+\varepsilon}
\]
uniformly for all $1/2\le\sigma\le 1$. These observations lead to the formulation of the \emph{Density Hypothesis} and, in turn, a new problem.
\begin{density}
For all $1/2\le \sigma \le 1$, we have $N(\sigma,T)\ll T^{2(1-\sigma)+\varepsilon}$ as $T\to \infty$.
\end{density}
 
\begin{problem}\label{problemB} Determine the least exponent $A(\sigma)$ for which one has the zero-density estimate 
$N(\sigma,T)\ll T^{A(\sigma)(1-\sigma)+\varepsilon}$ 
for all $1/2\le \sigma \le 1$. 
\end{problem}
\noindent We remark that in the notation of Problem \ref{problemB}, the Density Hypothesis asserts that $A(\sigma)\le 2$ for all $1/2\le \sigma\le 1$.

\subsection{Timeline of a few key zero-density theorems} We now highlight a few key results concerning Problem \ref{problemB}; we will discuss these results more thoroughly in the coming sections. In 1914, Bohr and Landau ~\cite{BL14} proved the first zero-density theorem, which states that for every \emph{fixed} $\sigma>1/2$,  $N(\sigma,T)\ll T$. While this bound is not sensitive to the value of $\sigma$, it does tell us that the zeros are necessarily clustered around the critical line. In particular, since $N(T)\sim (T/2\pi)\log T$,  the proportion of zeros of $\zeta(s)$ for which $1/2 < \beta <1$ tends to $0$ as $T\to \infty$. 

In 1940, Ingham ~\cite{Ingham40}, building off the work of Carlson~\cite{Carl21}, proved the following zero-density theorem. 
\begin{theorem}[Ingham, 1940]\label{thm: ingham1940}The estimate
\begin{equation}\label{eqn: ingham40}
N(\sigma,T)\ll T^{\frac{3(1-\sigma)}{2-\sigma}+\varepsilon}
\end{equation}
holds uniformly for $1/2 \le \sigma \le 1$ as $T\to \infty$. 
\end{theorem}

In 1969, Montgomery ~\cite{Montgomery69mean, Montgomery69} introduced a new approach to prove zero-density estimates, which recovers Ingham's bound \eqref{eqn: ingham40} in the range $1/2\le\sigma \le 4/5$, and  improves it in the range $4/5 < \sigma \le 1$. His work builds upon ideas of Halász \cite{Halasz68}, which was also used in \cite{HalTur69} to prove that the Lindelöf Hypothesis implies $N(\sigma,T)\ll T^\varepsilon$ for $3/4< \sigma \le 1$. (See also \cite[Theorem 2]{Montgomery69}.) Soon after, Huxley ~\cite{Huxley72} modified Montgomery's techniques, recovering Ingham's bound \eqref{eqn: ingham40} in the range $1/2\le\sigma \le 3/4$, and improving Montgomery's bound in the range $3/4 < \sigma \le 1$.
\begin{theorem}[Huxley, 1972]\label{thm: huxley}The estimate
\[
N(\sigma,T)\ll T^{\frac{3(1-\sigma)}{3\sigma-1}+\varepsilon}
\]
holds uniformly for $3/4\le \sigma \le 1$ as $T\to \infty$.
\end{theorem}
Taken together, Theorem \ref{thm: ingham1940} and Theorem \ref{thm: huxley} show that Problem \ref{problemB} holds for $A(\sigma)\le 12/5=2.4$ over the range $1/2\le \sigma \le 1$. That is, for all  $1/2\le \sigma \le 1$,
\begin{equation}\label{thm: huxrecord}
N(\sigma, T) \ll T^{\frac{12}{5}(1-\sigma)+\varepsilon}, \quad T\to \infty.
\end{equation}
This was the strongest bound known for all $1/2\le\sigma \le 1$ until very recently.
\begin{figure}[ht]
    \centering
    \includegraphics[width=0.8\linewidth]{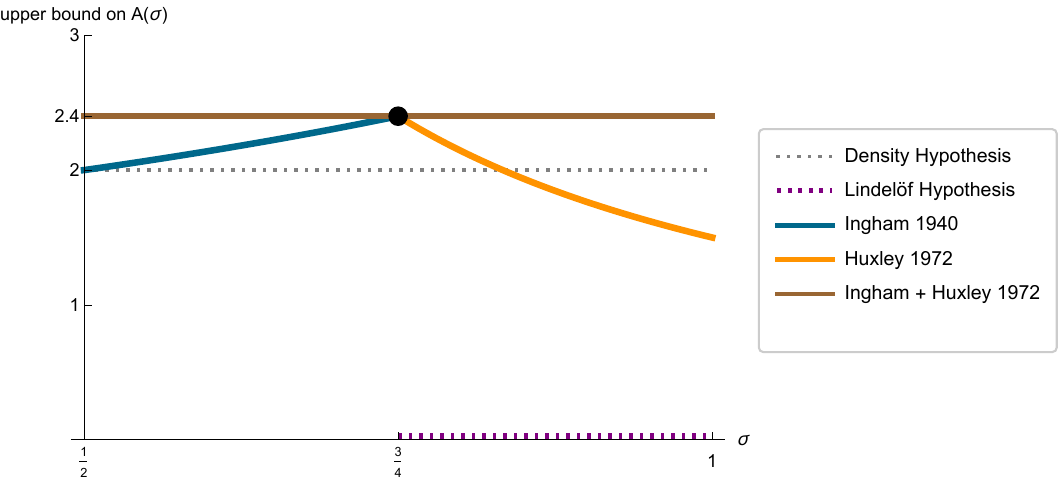}
    \caption{Upper bounds on $A(\sigma)$ for all $1/2\le \sigma \le 1$, prior to 2024}
\end{figure}

On May 31, 2024, Larry Guth and James Maynard  announced ~\cite{GuthMaynard24} a new zero-density theorem which gives the first improvement towards Problem \ref{problemB} in the range $1/2 \le \sigma \le 3/4$ since the 1940 estimate of Ingham \eqref{eqn: ingham40}. (That's an 84-year gap!) In particular, they prove the following, now published in \cite{GuthMaynard2026}.

\begin{theorem}[Guth-Maynard, 2026]\label{thm: guth-maynard}
The estimate
\[
N(\sigma,T)\ll T^{\frac{15(1-\sigma)}{3+5\sigma}+\varepsilon}
\]
holds uniformly for $1/2\le \sigma \le 1$ as $T\to \infty$.
\end{theorem}
Combining this estimate with Theorem \ref{thm: ingham1940} when $\sigma \le 7/10$, we see that Problem \ref{problemB} holds for $A(\sigma)<30/13\approx 2.308$ over the range $1/2\le\sigma \le 1$, that is 
\[
N(\sigma, T) \ll T^{\frac{30}{13}(1-\sigma)+\varepsilon}, \quad T\to \infty.
\]
This improves Huxley's record \eqref{thm: huxrecord} which resisted improvement for more than 50 years. 
\begin{figure}[ht]
    \centering
    \includegraphics[width=0.8\linewidth]{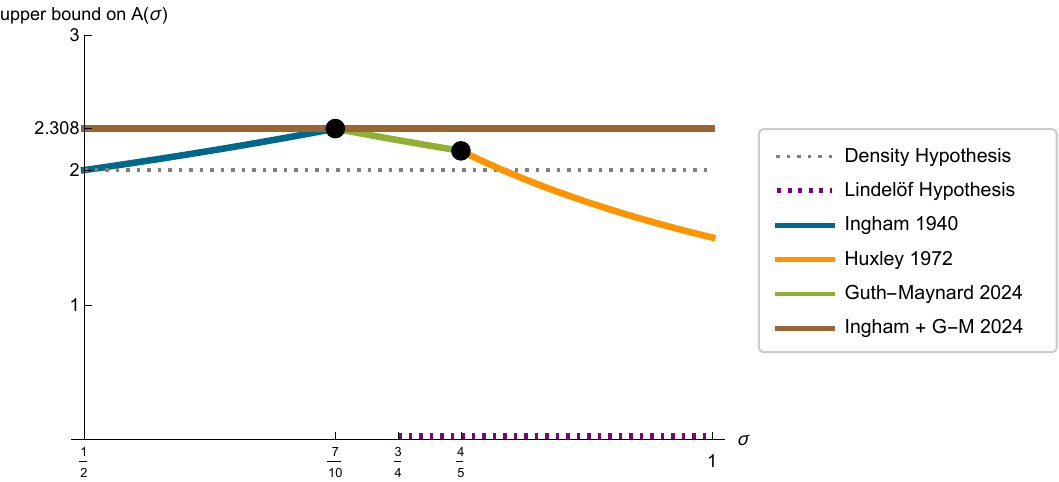}
    \caption{Upper bounds on $A(\sigma)$ for all $1/2\le \sigma \le 1$}
\end{figure}

\noindent\emph{Remark.} A detailed account of the strongest bounds on $A(\sigma)$ for various subintervals of $[1/2,1]$ may be found in ~\cite{GafniTao25}.

It is well-known (see, for example ~\cite{IK04} and Section 13 of ~\cite{GuthMaynard24}) that improvements on Problem \ref{problemB}  directly lead to improvements on Problem \ref{problemA}. In fact, improvements on Problem \ref{problemB} also yield results on the Prime Number Theorem in \emph{almost} all short intervals. We now state this relationship precisely, as in ~\cite[Theorem 1.1]{GafniTao25}.

\begin{theorem}\label{thm: correspondence}Let $A_0>0$ be such that for all $1/2 \le \sigma <1$, $A(\sigma)\le A_0$. Let $0 <\theta <1$ be fixed, and set $y=x^\theta$. 
\begin{enumerate}[{(i)}]
\item If $\theta > 1-1/A_0$, then \eqref{eqn: pntshort}  holds for all $x$.
\item If $\theta > 1-2/A_0$, then \eqref{eqn: pntshort}  holds for  all $x$ outside of an exceptional set of density zero.
\end{enumerate}
\end{theorem}

Theorem \ref{thm: correspondence} implies the following immediate corollary to Theorem \ref{thm: guth-maynard}.

\begin{corollary}[Guth-Maynard, 2026]
For large $x$, the Prime Number Theorem holds in  short intervals $(x,x+x^{17/30}]$ and in almost all intervals $(x,x+x^{2/15}]$.
\end{corollary}

\section{Zero-density estimates via Littlewood's Lemma} How does one go about proving a zero-density estimate, and what made Ingham's bound resist improvement for so long? Let's start with the first zero-density estimate obtained by Bohr and Landau \cite{BL14} in 1914 which states that for any fixed $\sigma>1/2$, $N(\sigma, T)\ll T$.

We begin with a special case of Littlewood's Lemma ~\cite[Theorem 1]{Littlewood24}  which, in analogy to Jensen's Formula (see, e.g., \cite[p.135]{SS}), concerns the zeros of a holomorphic function in a rectangle. 
\begin{littlewood}
Suppose that $f(s)$ is analytic and nonzero on the rectangle  with vertices $\sigma_0, \sigma_1, \sigma_1+ iT,$ and $\sigma_0+iT$, where $\sigma_0<\sigma_1$. Let $\nu(\sigma,T)$ denote the number of zeros of $f(s)$ in the part of the rectangle $\mathcal{C}$ for which $\re(s)> \sigma$. Then
\begin{align*}
2\pi \int_{\sigma_0}^{\sigma_1}\nu(\sigma,T)\,d\sigma &= \int_{0}^{T}\log |f(\sigma_0+it)|\,dt - \int_{0}^{T}\log |f(\sigma_1+it)|\,dt\\ &\hspace{.5in}+ \int_{\sigma_0}^{\sigma_1}\Arg f(\sigma+iT)\,d\sigma -\int_{\sigma_0}^{\sigma_1}\Arg f(\sigma)\,d\sigma.\notag
\end{align*}
\end{littlewood}

\subsection{Proof of the first zero-density theorem} We now sketch the proof of Bohr and Landau's zero-density theorem us Littlewood's Lemma. We follow the argument of \cite[Section 9.15]{Titchmarsh}.

Let $f(s)=\zeta(s)$, $T>0$, and $\sigma_1=2$. Fix $1/2< \sigma_0 \le 1$, where $\zeta(\sigma_0+it)\ne 0$ for $0\le t \le T$, and fix $\alpha=1/2+1/2(\sigma_0-1/2)$. Note that $1/2 < \alpha < \sigma_0 < 1$. We will show $N(\sigma_0,T)\ll T$.

Since $N(\sigma,T)$ is nonincreasing in $\sigma$, we have
\[
N(\sigma_0,T)\le \frac{1}{\sigma_0-a}\int_{\alpha}^{\sigma_0}N(\sigma,T)\,d\sigma \le \frac{2}{\sigma_0-1/2}\int_{\alpha}^{1}N(\sigma,T)\,d\sigma.
\]
By Littlewood's Lemma and the fact that $N(\sigma,T)=0$ for $\sigma>1$, we have
\begin{equation}\label{eqn}
\int_{\alpha}^{1}N(\sigma,T)\,d\sigma =\frac{1}{2\pi}\int_{0}^{T}\log |\zeta(\sigma_0+it)|\,dt +O(\log T).
\end{equation}
(For details on obtaining the error term, see ~\cite[p.229]{Titchmarsh}). To obtain a suitable upper bound on the integral on the right-hand side of \eqref{eqn}, we apply Jensen's Inequality (see, e.g., \cite[p.115]{Royden}) to find
\[
\frac{1}{2\pi}\int_{0}^{T}\log |\zeta(\sigma_0+it)|\,dt = \frac{T}{4\pi}\left(\frac{1}{T}\int_{0}^{T}\log \left(|\zeta(\sigma_0+it)|^2\right)\,dt\right) \le \frac{T}{4\pi}\log\left(\frac{1}{T} \int_{0}^{T}|\zeta(\sigma_0+it)|^2\,dt\right).
\]
Thus
\[
N(\sigma_0,T) \le \frac{T}{2\pi(\sigma_0-1/2)}\log\left(\frac{1}{T} \int_{0}^{T}|\zeta(\sigma_0+it)|^2\,dt\right) + O(\log T).
\]
We have reduced the problem to another central topic in analytic number theory: mean values of the Riemann zeta-function. Since $\sigma_0>1/2$, Lemma 1 of \cite{Littlewood24} yields
\begin{equation}\label{eqn: certain}
\int_{0}^{T}|\zeta(\sigma_0+it)|^2\,dt \ll  \frac{A}{\sigma_0-1/2}T,
\end{equation}
where $A>0$ is an absolute constant. Assembling these results, we find  
\[
N(\sigma_0, T) \ll T,
\]
where the implied constant  depends on $\sigma_0$.

\subsubsection{Historical note} The argument given above is anachronistic: Bohr and Landau's result was published in 1914, nearly 10 years before Littlewood's Lemma appeared. (The proof above, in fact, appears in \cite{Littlewood24}.) Bohr and Landau proved their result using Jensen's Formula, which we now state for comparison. Let $f(s)$ be an analytic function in a region of the complex plane which contains the closed disc $D_R$ of radius $R>0$ about the origin. Suppose that $f(0)\ne0$ and $f(s)$ is nonzero on the boundary of $D_R$. Let $n(t)$ denote the number of zeros of $f(s)$ in the disc of radius $t$ centered at the origin. Then
\[
\int_{0}^{R}\frac{n(r)}{r}\, dr=\frac{1}{2\pi}\int_{0}^{2\pi}\log|f(Re^{i\theta})|\,d\theta -\log|f(0)|.
\]
As discussed in ~\cite[p.250--251]{IK04}, to count the number of zeros of $f(s)$ in $D_R$, it is helpful to apply Jensen's Formula not to $f(s)$ directly, but instead to a product $f(s)g(s)$, where $g(s)$ is holomorphic in $D_R$ and large when $f(s)$ is small. In the case of $f(s)=\zeta(s)$, Bohr and Landau took $g(s)$ to be a truncation of the reciprocal of the Euler product \eqref{def: eulerprod} of $\zeta(s)$, i.e. $g(s)=\prod_{p\le X}(1-p^{-s})$. As we will see in the next sections, subsequent improvements on upper bounds on $N(\sigma,T)$ employ a similar idea, replacing the truncated Euler product with its Dirichlet series counterpart. 

\subsection{Ingham's 1940 zero-density theorem} Ingham also made use of Littlewood's Lemma to prove his zero-density theorem, so we should expect some sort of mean value theorem involving $\zeta(s)$ to appear in the proof of Theorem \ref{thm: ingham1940}. We focus on this aspect of Ingham's result in this section. As Ingham writes in an earlier paper ~\cite{Ingham37} on the subject, ``The standard proofs of theorems of the type [in Problem \ref{problemB}] are based on the estimation of a certain integral\ldots.'' This integral is similar to that of \eqref{eqn: certain}, with $\zeta(s)$ replaced by a related function constructed as follows. 

Let $\mu$ denote the Möbius function, where $\mu(1):=1$ and for $n>1$,
\[
\mu(n):=\begin{cases}
(-1)^k & \text{if } n \text{ is the product of } k \text{ distinct primes},\\
0 & \text{if } n \text{ is divisible by a square greater than }1.
\end{cases}
\]
We will make use of the following two important properties of the Möbius function in this section:
\[
\sum_{d|n}\mu(d) = \begin{cases}
1 & \text{if } n=1,\\
0 & \text{if } n>1,
\end{cases}
\qquad \text{and} \qquad
\frac{1}{\zeta(s)}=\sum_{n=1}^{\infty}\frac{\mu(n)}{n^s},\qquad \re(s)>1.
\]
Let $X$ be large, and define
\[
M_X(s):=\sum_{m\le X}\frac{\mu(m)}{m^s}.
\]
For $\re(s)>1$,
\begin{equation}\label{def: zetaMX}
\zeta(s)M_X(s) = \sum_{k=1}^{\infty}\frac{1}{k^s}\sum_{m\le X}\frac{\mu(m)}{m^s}= \sum_{k=1}^{\infty}\sum_{m\le X}\frac{\mu(m)}{(mk)^s} = \sum_{n=1}^{\infty}\frac{a_n}{n^s}, \quad \text{where }\, a_n:=\sum_{\substack{m|n\\m\le X}}\mu(m).
\end{equation}
Note that $a_1=1$ and $a_n=0$ for $2 \le n \le X$, and observe that if $\zeta(\rho)=0$ then $\zeta(\rho)M_X(\rho)=0$ too, and that for large $X$, we might expect that $\zeta(s)M_X(s)\approx 1$ for $1/2\le \re(s) \le 1$. 

Define
\[
f_X(s):=\zeta(s)M_X(s)-1.
\]
The ``certain integral" to which Ingham refers is
\begin{equation}\label{eqn: certainintegral}
\int_{0}^{T}|f_X(\sigma+it)|^2\,dt.
\end{equation}
The need for an upper bound on \eqref{eqn: certainintegral} arises in the course of applying Littlewood's Lemma to the function $h_X(s):=1-(f_X(s))^2$, since $\log|h_X(s)| \le \log(1+|f_X(s)|^2) \le |f_X(s)|^2$. In particular, since $h_X(s)$ can be written as the product of $\zeta(s)$ and $M_X(s)(2-\zeta(s)M_X(s))$, a zero-density estimate for $h_X(s)$ yields an upper bound for that of $\zeta(s)$.

Ingham proves an upper bound on \eqref{eqn: certainintegral} that holds uniformly on $1/2\le \sigma \le 1$  via a convexity theorem for mean values of analytic functions. In practice this involves finding a suitable upper bound on \eqref{eqn: certainintegral} with $\sigma=1/2$ and another on \eqref{eqn: certainintegral} with $\sigma=1+\delta$, where $\delta>0$. His argument leads to the estimate $N(\sigma,T)\ll T^{(1+2\sigma)(1-\sigma)}\log^5 T$, which holds uniformly for $1/2 \le \sigma \le 1$ as $T\to \infty$. (See  \cite[Theorem 5]{Ingham37}.)  A few years later, Ingham ~\cite{Ingham40} improved this result by proving a more flexible version of \eqref{eqn: certainintegral}.

\begin{theorem}[Ingham, 1940]\label{thm: inghammeanvalue}
Let $T>1, X\ge 3,$ and $1/2\le \sigma \le 1$. Then
\[
\int_{1}^{T}|f_X(\sigma+it)|^{u}\,dt \ll (T+X)X^{-(2\sigma-1)/(2-\sigma)}\log^4(T+X),
\]
for a certain\footnote{A precise description of $u$ can be found on page 292 on \cite{Ingham40}.} $u=u(\sigma, T, X)$ in the range $4/3 \le u \le 2$.
\end{theorem}
Ingham proves this bound using a two-variable convexity theorem of Gabriel ~\cite{Gabriel1927} and interpolates between an upper bound on $\int_{0}^{T}|f_X(1+\delta+it)|^{2}\,dt$, where $ 0 < \delta < 1$, and an upper bound on
\begin{equation}\label{eqn: inghammeanvalue}
\int_{0}^{T}|f_X(1/2+it)|^{4/3}\,dt.
\end{equation}

We briefly sketch how to bound \eqref{eqn: inghammeanvalue}, as it is a key piece of Ingham's overall argument. (It will also quickly lead to the topic of mean values of Dirichlet polynomials, which we will later explore in the context of Huxley's work.) To begin, Ingham uses the triangle inequality and Hölder's inequality with $p=3$ and $q=3/2$ to show
\begin{align*}
\int_{0}^{T}|f_X(1/2+it)|^{4/3}\,dt &\ll \int_{0}^{T}\left(|\zeta(1/2+it)|^{4/3}|M_X(1/2+it)|^{4/3}+1\right)\,dt\\
&\hspace{.5in}\ll T+\left(\int_{0}^{T}|\zeta(1/2+it)|^4\,dt \right)^{1/3}\left(\int_{0}^{T}|M_X(1/2+it)|^2\,dt \right)^{2/3}.
\end{align*}
Hardy and Littlewood ~\cite{HardyLittlewood1923} showed that $\int_{0}^{T}|\zeta(1/2+it)|^4\,dt \ll T\log^4 T$. To bound the remaining integral, we appeal to the following classical mean value theorem. (See \cite[Theorem 9.1]{IK04}.)
\begin{theorem}[Mean Value Theorem for Dirichlet Polynomials]\label{thm: imvt}
For any complex numbers $b_n$, we have
\[
\int_{0}^{T}\left|\sum_{1\le n \le N}b_nn^{it} \right|^2\,dt = (T+O(N))\sum_{1\le n \le N}|b_n|^2,
\]
where the implied constant is absolute.
\end{theorem}
Thus
\[
\int_{0}^{T}|M_X(1/2+it)|^2\,dt\ll (T+X)\sum_{1\le n \le X}\frac{1}{n}\ll (T+X)\log X,
\]
and hence 
\begin{align*}
\int_{0}^{T}|f_X(1/2+it)|^{4/3}\,dt &\ll (T+X)\log^2(T+X).
\end{align*}
For full details of Ingham's proof of Theorem \ref{thm: inghammeanvalue}, see ~\cite{Ingham40} and ~\cite{Ingham37}.

\section{Zero-density estimates via zero-detecting polynomials}
Twenty-five years after Ingham's paper, A.~I. Vinogradov \cite{Vinogradov1965} introduced a new method to estimate $N(\sigma,T)$, which Montgomery \cite{Montgomery69} developed just a few years later. Unlike the arguments presented above, the method does not approach counting zeros via Littlewood's Lemma. Instead, it counts zeros by detecting how often a particular Dirichlet polynomial is large. We will see that the method produces much stronger bounds when treating zeros close to the line $\re(s)=1$. 

We now set up the zero-detecting method, combining aspects of the expositions of \cite{Montgomery69} and \cite{MaynardPratt24}. Let $T$ be large, and let $X\le Y \le T^A$  for some absolute constant $A>0$. We will specify $X$ and $Y$ in terms of $T$ at the end of the argument. Consider the Dirichlet series
\begin{equation}\label{eqn: IZ}
  I(s):=\sum_{n=1}^{\infty}\frac{a_n}{n^s}e^{-n/Y}, 
\end{equation}
where $a_n$ is as in \eqref{def: zetaMX}. As previously noted,  $a_1=1$ and $a_n=0$ for $2\le n\le X$. Thus
\begin{equation}\label{eqn: IZ1}
  I(s) = e^{-1/Y}+\sum_{n>X}\frac{a_n}{n^s}e^{-n/Y}. 
\end{equation}

We may express $I(s)$ another way using the inverse Mellin transform of the Gamma function, which for $\re(z)>0$ and $c>0$ states
\[
e^{-z}= \frac{1}{2\pi i}\int_{(c)}\Gamma(w)z^{-w}\,dw. 
\]
Applying this to \eqref{eqn: IZ} and taking $c=2$, say, we find
\begin{equation}\label{eqn: mellin}
  I(s) =\frac{1}{2\pi i}\int_{(2)}\zeta(s+w)M_X(s+w)Y^w\Gamma(w)\,dw.
\end{equation}
Let $\rho=\beta+i\gamma$ denote a zero of $\zeta(s)$ counted by $N(\sigma,T)$. We may assume that $\beta \ge \sigma \ge 1/2 + 1/\log T$, and we recall that $0<\gamma \le T$. Evaluate $I(\rho)$ in \eqref{eqn: mellin} and then shift the contour to $\re(w)=1/2-\beta$. We pass through a pole at $w=1-\rho$ from $\zeta(\rho+w)$ with residue $M_X(1)Y^{1-\rho}\Gamma(1-\rho)$, while the pole of $\Gamma(w)$ at $w=0$ is canceled by the zero of $\zeta(\rho+w)$ there. Comparing the resulting expression to \eqref{eqn: IZ1} evaluated at $s=\rho$, we find
\begin{align*}
e^{-1/Y}+\sum_{n>X}\frac{a_n}{n^\rho}e^{-n/Y} = M_X(1)Y^{1-\rho}\Gamma(1-\rho) 
+\frac{1}{2\pi i}\int_{(\frac{1}{2}-\beta)}\zeta(\rho+w)M_X(\rho+w)Y^{w}\Gamma(w)\,dw + o(1).
\end{align*}
Consider the left-hand side of the above equation. Since $Y$ is large, $e^{-1/Y}=1+o(1)$, and the tail of the Dirichlet series is negligible since the factor $e^{-n/Y}$ decays rapidly. On the right-hand side, note that if $\gamma\ge (\log T)^2$ then $M_X(1)Y^{1-\rho}\Gamma(1-\rho)=o(1)$ due to the rapid decay of the Gamma function in vertical strips. Therefore, for any $\rho$ for which $\gamma\ge (\log T)^2$ we have
\begin{equation}\label{eqn: key}
1+\sum_{X < n \le Y^2}\frac{a_n}{n^\rho}e^{-n/Y} =
 \frac{1}{2\pi i}\int_{(\frac{1}{2}-\beta)}\zeta(\rho+w)M_X(\rho+w)Y^{w}\Gamma(w)\,dw +o(1).
\end{equation}

We are now able to make a key observation: because of the 1 appearing on the left-hand side of \eqref{eqn: key}, it is not possible for the Dirichlet polynomial and the integral to be simultaneously small. Splitting up the Dirichlet polynomial into $\ll \log T$ dyadic sub-intervals, we can rephrase this observation to state that, apart from the zeros of $\zeta(s)$ with $0<\gamma\le (\log T)^2$, either

\begin{equation}\label{eqn: type1}
\left | \sum_{N< n \le 2N}\frac{a_n}{n^\rho}e^{-n/Y}\right | \ge \frac{1}{3\log T} \qquad \text{for some } N \in \left(X\,,\,Y^2\right]
\end{equation}
or
\begin{equation}\label{eqn: type2}
\left|\frac{1}{2\pi i}\int_{(1/2-\beta)}\zeta(\rho+w)M_X(\rho+w)Y^w\Gamma(w)\,dw\right| \ge \frac{1}{3}.
\end{equation}
A zero that satisfies \eqref{eqn: type1} is called a Type~I zero, and a zero that satisfies \eqref{eqn: type2} is called a Type~II zero. Let $\mathcal{R}$ denote the subset of the zeros counted by $N(\sigma,T)$ for which $|\gamma_1-\gamma_2|\ge 1,\,\rho_1,\rho_2\in \mathcal{R}$. Such a set is called \emph{well-spaced} or \emph{1-separated}. By \eqref{eqn: NT}, there are $O(\log T)$ nontrivial zeros with $\gamma\in (t,t+1]$ for any $t \in [0, T]$, so $\mathcal{R}$ is nonempty. Let $R_1, R_2$ denote the number of Type~I, Type~II zeros in $\mathcal{R}$, respectively. Since there are $\ll \log^3T$ zeros counted by $N(\sigma,T)$ for which $0<\gamma\le (\log T)^2$, we have
\[
N(\sigma,T)\ll R_1(\log T)^2 + R_2(\log T)^2+(\log T)^3.
\]
The contribution of the Type~II zeros can be controlled  (see, e.g. \cite[Lemma 6.3]{MaynardPratt24}), while that of the Type~I zeros is more delicate. Thus, towards bounding $R_1$, we now aim to understand the frequency for which Dirichlet polynomials over well-spaced sets can be large. 

\section{Large values of Dirichlet polynomials}
Let $N\ge 1$, $(b_n)_{N< n \le 2N}$ a sequence of complex numbers, and consider the Dirichlet polynomial
\[
D_N(t):=\sum_{N\le n \le 2N}b_nn^{it}.
\]
Let $W$ be a set of 1-separated points in $[0,T]$ with $T\ge 1$, and suppose that there exists $V>0$ such that $|D_N(t)|\ge V$ for every $t\in W$. How big can $|W|$ be under these circumstances? By definition, $W$ contains at most $T+1$ points. On the other hand, 
\begin{equation*}
|W| = \sum_{t\in W}1 \le V^{-2}\sum_{t\in W}|D_N(t)|^2.
\end{equation*}
Thus we can obtain a nontrivial bound on $|W|$ via a mean value theorem for $D_N(t)$ over $W$. In 1967, Gallagher \cite{Gallagher67} proved a discrete version of Theorem \ref{thm: imvt}, which in the context of Dirichlet polynomials may be stated as follows. (See \cite[Theorem 9.4]{IK04}.)
\begin{theorem}[{Discrete Mean Value Theorem}] \label{thm: discretemean} Let $W$ be a set of 1-separated points in $[0,T]$ with $T\ge 1$, and let $(b_n)$ be any sequence of complex numbers. Then we have
\[
\sum_{t\in W}\Big|\sum_{N\le n \le 2N}b_nn^{it}\Big|^2 \ll (T+N)(\log 2N)\sum_{N\le n \le 2N}\left|b_n\right|^2,
\]
where the implied constant is absolute.
\end{theorem}

\subsection{Recovering Ingham's zero-density estimate} We now show how to recover Ingham's zero-density estimate using large values of Dirichlet polynomials. In particular, we focus on obtaining an upper bound on $R_1$, the number of 1-separated Type~I zeros in $\mathcal{R}$ counted by $N(\sigma,T)$. We follow the argument of Montgomery \cite{Montgomery69mean}. 

Set $b_n=a_nn^{-\beta}e^{-n/Y}$. By \eqref{eqn: type1} and Theorem~\ref{thm: discretemean}, there exists an $N\in(X,Y^2]$ such that
\begin{equation}\label{eqn: r1bound}
R_1  \ll (T+N)(\log 2N)(\log T)^2\sum_{N\le n \le 2N}\left|b_n\right|^2.
\end{equation}
By the definition of $a_n$ given in \eqref{def: zetaMX}, $|a_n| \le d(n)$, where $d(n)$ denotes the divisor function which counts the number of divisors of an integer $n\ge 1$. 
It is a result of Dirichlet (see, e.g., \cite[Theorem 3.3]{Apostol76}) that for all $x\ge 1$, 
\[
\sum_{n\le x}d(n)=x\log x +O(x).
\]
By induction and partial summation one can show that for integers $k\ge 1$, 
\begin{equation}\label{eqn: divisorsum}
\sum_{n\le x}d(n)^k\ll x(\log x)^{2^k-1}.
\end{equation}
(See \cite[Theorem 5.3]{Hua}).) Recall that $\beta \ge \sigma$ and $X\le N \le Y^2$. Applying \eqref{eqn: divisorsum} with $k=2$,  we find
\begin{align*}
\sum_{N\le n \le 2N}\left|b_n\right|^2 &\ll \sum_{N\le n \le 2N}|a_n|^2n^{-2\beta}e^{-2n/Y} \ll N^{-2\sigma}e^{-2X/Y}N(\log 2N)^3.
\end{align*}
Inserting this bound into \eqref{eqn: r1bound} yields
\[
R_1 \ll  (TX^{1-2\sigma}+Y^{2-2\sigma})(\log T)^{O(1)}.
\]

To deduce an upper bound on $R_2$, we further partition the well-spaced Type~II zeros, i.e., the zeros satisfying \eqref{eqn: type2}, depending on the size of $|M(\rho+w)|$. Using Theorem~\ref{thm: discretemean} and the upper bound on the fourth moment on $\zeta(1/2+it)$ referenced above the statement of Theorem~\ref{thm: imvt}, one can show
\[
R_2 \ll TY^{\frac{2-4\sigma}{3}}(\log T)^{O(1)}.
\]
(See \cite[p. 352]{Montgomery69} for the details on obtaining this bound.) Combining these estimates, we find
\[
N(\sigma,T) \ll  \left(TX^{1-2\sigma}+Y^{2-2\sigma}+ TY^{\frac{2-4\sigma}{3}}\right)(\log T)^{O(1)}.
\]
Take $X=T$ and $Y=T^{\frac{3}{2(2-\sigma)}}$, which balances the two terms involving $Y$. We have
\[
N(\sigma,T) \ll  T^{\frac{3(1-\sigma)}{2-\sigma}}(\log T)^{O(1)},
\]
which recovers Ingham's bound \eqref{eqn: ingham40} up to logarithms. \\

\subsection{Montgomery-Halász-Huxley large values estimate}
In 1969, Montgomery ~\cite{Montgomery69mean} used a method of Halász ~\cite{Halasz68} to prove a new mean value theorem for Dirichlet polynomials over a set $W$ of 1-separated points. Using the duality principle of bilinear forms and bounds on exponential sums, Montgomery's mean value theorem produces an upper bound that, in contrast to Theorem \ref{thm: discretemean}, is sensitive to the cardinality of $W$. In particular, the bound is much stronger when $|W|$ is smaller than $T^{1/2}$.

\begin{theorem}[Montgomery-Halász Mean Value Theorem] Suppose $(b_n)$ is a sequence of complex numbers and  $(t_r)_{r\le R}$ is a sequence of 1-separated points in $[0,T]$ with $T\ge 1$. Then we have
\[
\sum_{t_r}\Big|\sum_{N\le n \le 2N}b_nn^{it_r}\Big|^2 \ll (N+RT^{1/2})\log(2T)\sum_{N\le n \le 2N}\left|b_n\right|^2,
\]
where the implied constant is absolute.
\end{theorem}

Refining the ideas in ~\cite{Montgomery69mean} and \cite{Montgomery69}, Huxley~\cite{Huxley72} proved what we now refer to as the Montgomery-Halász-Huxley large values estimate. 

\begin{theorem}[Montgomery-Halász-Huxley large values estimate]\label{thm: MHH} Suppose $(b_n)$ is a sequence of complex numbers with $|b_n|\le 1$ and $(t_r)_{r\le R}$ is a sequence of 1-separated points in $[0,T]$ such that
\[
\Big|\sum_{N\le n \le 2N}b_nn^{it_r}\Big|\ge V
\]
for some $V>0$. 
Then 
\begin{equation*}
R\le T^{o(1)}\left(N^2V^{-2}+T\min\left\{NV^{-2},N^4V^{-6}\right\}\right),
\end{equation*}
where the implied constant is absolute.
\end{theorem}
We remark that Theorem \ref{thm: MHH} yields much better bounds on Type~I zeros for which $\beta \ge \sigma \ge 3/4$ but cannot do so for $1/2\le \sigma <3/4$ due to the framework on which it is built.

\subsection{Guth-Maynard large values estimate}
In 2024, Guth and Maynard \cite{GuthMaynard24} introduced a new method to obtain large value estimates of Dirichlet polynomials. Their remarkable work gives a new way to count Type~I zeros, which in turn provides a stronger zero-density theorem for $7/10 \le \sigma \le 8/10$. It is worth stating again that this gives the first progress for $\sigma \le 3/4$ since 1940! 

\begin{theorem}[Guth-Maynard large values estimate]\label{thm: newbound} Suppose $(b_n)$ is a sequence of complex numbers with $|b_n|\le 1$ and $(t_r)_{r\le R}$ is a sequence of 1-separated points in $[0,T]$ such that
\[
\Big|\sum_{N\le n \le 2N}b_nn^{it_r}\Big|\ge V
\]
for some $V>0$. Then 
\[
R\le T^{o(1)}\left(N^2V^{-2}+N^{18/5}V^{-4}+TN^{12/5}V^{-4}\right),
\]
where the implied constant is absolute.
\end{theorem}
\subsubsection{Obtaining a new zero-density theorem} As Guth and Maynard explain \cite[page 45]{GuthMaynard24}, Theorem \ref{thm: newbound} is applied when the length of the relevant Dirichlet polynomial satisfies
\[
T^{5/(3+5\sigma)}\le N^2 \le T^{75(1-\sigma)/(54+30\sigma-100\sigma^2)}.
\]
The term $N^2V^{-2}$ is dominated by $N^{18/5}V^{-4}$ in the critical range $7/10\le \sigma \le 8/10$, and in this regime they show
\[
R_1 \ll T^{15(1-\sigma)/(3+5\sigma)}.
\]
Finally, for $N^2$ outside of the range above, the authors use classical estimates to bound $R_1$. Combining these results, the authors ultimately obtain the zero-density estimate
\[
N(\sigma,T) \ll T^{15(1-\sigma)/(3+5\sigma)+o(1)}.
\]

\subsubsection{The additive energy of large value sets} The authors observe (see \cite[page 4]{GuthMaynard24}) that one might expect to have improved bounds for Dirichlet polynomials whose large value set $W$ has a lot structure.  
For example, if the points of $W$ form an arithmetic progression, it seems reasonable that this information could be harnessed to obtain improved bounds. They make this precise in terms of the \emph{additive energy} of a finite set $W$, which is defined as
\[
E(W):=\#\{w_1, w_2, w_3, w_4 \in W\,:\, |w_1+w_2-w_3-w_4|< 1\}.
\]
They show that a result of Heath-Brown \cite[Theorem 1]{HB79} implies that one may indeed obtain stronger bounds for Dirichlet polynomials whose large value set $W$ has ample additive energy. For sets with low additive energy, Guth and Maynard give a new method to obtain good bounds.

We now follow Section 2 of \cite{GuthMaynard24} to sketch the set-up of the method and highlight a few of their key remarks. Let $N\ge 1$ denote the length of the Dirichlet polynomial $D_N(t)$ and $W$ a set of 1-separated points. As we have seen above, we can obtain an upper bound on $|W|$ by finding an upper bound on 
\[
\sum_{t\in W}|D_N(t)|^2.
\]
We can rephrase this task in terms of bounding singular values of matrices. Let $M$ denote a $|W|\times N$ matrix for which $M_{t,n}= n^{it}$. Then for $\vec{b}\in \mathbb{R}^N$, 
\[
(M\vec{b})_t=\sum_{N\le n \le 2N}b_nn^{it} = D_N(t).
\]
Thus 
\[
\sum_{t\in W}|D_N(t)|^2 = \lVert M\vec{b}\rVert^2 \le N\lVert M\rVert^2, 
\]
where $\lVert M\rVert$ denotes the matrix norm of $M$. In particular,  $\lVert M\rVert = s_1(M)$ is the largest singular value of $M$. The authors describe how various reasonable approaches to bounding $s_1(M)$ would eventually fail, and then they state that one can prove a bound similar to 
\[
s_1(M)^6 \ll \Big| \, \sum_{\substack{t_1,t_2,t_3\in W \\ |t_j-t_k|>T^\epsilon\, \forall j,k}}\sum_{n_1,n_2,n_3\sim N}n_1^{i(t_1-t_2)}n_2^{i(t_2-t_3)}n_3^{i(t_3-t_1)}\,\Big|.
\]
The aim, then, is to show that the right-hand side is small. As the authors discuss, after applying Poisson summation above, the next maneuver would typically be to apply stationary phase. By refraining from taking this natural next step, the authors obtain an expression in which the variables $t_1, t_2, t_3$ above are separated from one another. In doing so, however, they have given up a factor of $T^{3/2}$. Nevertheless,  the authors show how to harness the additive structure of $W$, whether it be high or low, to obtain enough cancellation to win. \\

\noindent\textbf{Acknowledgments.}\,The author is partially supported by NSF CAREER DMS-2239681 and thanks Brian Conrey, Dan Goldston, Micah Milinovich, and the  referee for helpful comments and suggestions. 

\bibliographystyle{amsalpha}
\bibliography{references}
\end{document}